\documentclass[12pt]{article}

\usepackage[a4paper,left=2.5cm,top=3cm,textwidth=16cm,textheight=24.5cm]{geometry}

\usepackage[bookmarksopen=true,bookmarks=true,unicode,setpagesize]{hyperref}
\hypersetup{colorlinks=true,linkcolor=black,citecolor=black}

\usepackage[intlimits]{amsmath}

\usepackage{amsthm,amssymb}

\theoremstyle{plain}
\newtheorem{theorem}{Theorem}

\theoremstyle{definition}

\theoremstyle{remark}

\allowdisplaybreaks[3]

\renewcommand\C{\mathbb{C}}
\newcommand\M{\mathcal{M}}
\newcommand\N{\mathbb{N}}
\newcommand\R{\mathbb{R}}

\newcommand\Ga{\Gamma}
\newcommand\ga{\gamma}
\newcommand\la{\lambda}
\newcommand{\om}{\omega}

\newcommand\X{{\R^d}}

\begin{document}

\pagestyle{plain}
\title{Philosophy  of Natural Numbers}

\date{}
\author{
\textbf{Yuri Kondratiev}\\
 Department of Mathematics, University of Bielefeld, \\
 D-33615 Bielefeld, Germany,\\
 Dragomanov University, Kyiv, Ukraine\\
 Email: kondrat@math.uni-bielefeld.de
 }

\maketitle

\begin{center}
\begin{abstract}
We discuss an extension of classical combinatorics theory to the case
of spatially distributed objects. 

\end{abstract}



\end{center}

%\tableofcontents{}

\section{Introduction}

The set of natural numbers $\N=\{0,1,2, \dots\}$ is a fundamental object in the mathematics.
In certain sense $\N$ is the root of all modern mathematics. Other mathematical structures 
may be created as a logical development of this object. The latter motivated
L. Kronecker who summarized "God made the integers, all else is the work of man".There is famous citation from 
I.Kant: "Two things fill the mind: the starry heavens above me and the moral law within me".  A mathematician may continue: "and natural numbers given to my  mind".

From the time of Pithagoras philosophers was  trying to see hidden meaning of natural numbers and their mystical properties. Considering  $\N$ as a set of real things in
mathematics we will ask ourself about possible ideas behind these numbers.
The myth of Plato's Cave served as one of the motivations for developing his concept of a world of ideas and a world of things. In the dialogue "State" he gives several examples illustrating this concept. As we know, Plato considered mathematics as one of the most important building blocks used to construct his philosophical system. Mathematical theories can serve as simple and illustrative tools for the existence of a world of ideas and a world of things. In a number of model situations, we are dealing with objects that appear from our observations in physics, biology, ecology, etc., yet full understanding of the mathematical structures of these models requires consideration of more general mathematical theories, which under some canonical mapping leads to the model situations in question.

The first and essentially obvious observation here is the following. A number $n\in\N$
we interpret as a number of objects (a population)  located in a location space $X$. For simplicity
we take $X=\X$. The collection of all $n$-point subsets (or configurations with
$n$ elements) 
form a locally compact space $\Ga^{(n)}(\X)$. It is the space (quite huge)  of ideas 
for the number $n$.  Then to $\N$ corresponds the set
$$
\Ga_0(\X)= \cup_{n=0}^\infty \Ga^{(n)}(\X)
$$
of all finite configurations.
We can consider additionally the set $\Ga(\X)$ consisting all locally finite configurations.
This set may be considered as the space of ideas which corresponds to natural numbers
and additionally   to the actual infinity 
which is absent in the classical framework on natural numbers.

In such extension  of $\N$ we arrive in the main question. Namely, most important mathematical theories
related to natural numbers we need to develop to this new level. It concerns, first of all,  the combinatorics
that play central role in many mathematical structures and applications from probability theory to
genetics. In this note we will try to show such possibility trying to be as much as possible on
technically simple ground. To be friendly to more wide audience,
we restrict out explanations to descriptions of main constructions and 
formulation of some particular  results. For detailed discussions and extended references we refer to the recent
paper \cite{FKLO}.

\section{Classical combinatorics}

The combinatoric  is dealing with the set of natural numbers $\N$ and relations between them. As an important object we introduce binomial coefficients:
$$
\binom nk = \frac{n(n-1) \dots (n-k+1)}{k!}
$$ 
defined for $n\in\N$  and $0\leq   k\leq n$. Introducing the falling factorial  $(n)_k$ we can write 
$$
\binom nk= \frac {(n)_k}{k!}.
$$
These coefficients may be extended using embedding $\N\subset \R$ to polynomials
$$
N_k(t):= \binom tk = \frac{t(t-1)\dots (t-k+1)}{k!}= \frac{(t)_k}{k!},\; t\in\R
$$
which are called Newton polynomials.
For Newton polynomials hold Chu-Vandermond relations:
$$
(t+s)_n= \sum_{k=0}^n \binom nk (t)_k (s)_{n-k}.
$$

An alternative definition is given by the  generation function
$$
e_\la(t):=e^{t\log(1+\la)} =\sum_{n=0}^\infty \frac{\la^n}{n!} (t)_n= \sum_{n=)}^\infty \la^n N_n(t).
$$
Such transition to continuous variables makes possible to apply in discrete mathematics methods of
analysis. Note that 
using many particular generation functions we may create different polynomial systems.

Transition to continuous variables makes possible to apply in discrete mathematics methods of
analysis. In particular, let us define for functions $f:\R\to\R$ difference operators
$$
(D^+f)(t)= f(t+1)- f(t),
$$
$$
( D^- f)(t)= f(t-1) -f(t).
$$
By a direct computation we obtain
$$
D^+ (t)_n= n \; (t)_{n-1},
$$
$$
D^- (t)_n= -n \; (t-1)_{n-1}.
$$ 
Additionally,
$$
D^+ e_\la(t) = \la e_\la(t).
$$
In this way we arrive in the  framework of difference calculus closely related with the combinatorics \cite{DC}.
There are specific questions inside of difference calculus as, e.g., an analysis of Newton series
$$
\sum_{n=0}^\infty a_n N_n(t)
$$ and many others.

For functions $a:\N\to \R$ we define $b:\N\to\R$ as
$$
b=Ka,\;\; b(n)= \sum_{k=0}^n \binom nk a(k).
$$
This operator $K$  (aka combinatorial transform) is very useful in combinatorics and its inverse gives so-called
inclusion-exclusion formula:

$$
a(n)= \sum_{k=0}^n  \binom nk (-1)^{n-k} b(k).
$$

Note that for $a:\N\to \R,\;\; a(j)=0, j\neq k, a(k)=1$

$$
(Ka)(n)= \binom nk = k! N_k(n).
$$

\section{Spatial combinatorics}

Any $n\in \N$ we interpret as the size of a population. It is convenient in the study of
population models. There is a natural generalization leading to spatial ecological 
models. Now we would like to consider 
objects located in a given locally compact space $X$. For simplicity we will work with the
Euclidean space $\X$.  For the substitution of $\N$ in this situation we can use two possible sets.
Denote $\Ga(\X)$ the set of all locally finite configurations (subsets) from $\X$. 

$$
\Ga(\X) =\{ \ga\subset \X\;|\; |\ga\cap K|<\infty, {\rm{any\;\; compact} } \;\; K\subset \X\}.
$$

It is  the first version 
of the space in the spatial (continuous)  combinatoric we will use.

Another possibility, is to introduce the set
of all finite configurations $\Ga_0(\X)$.  Then
$$
\Ga_0(\X) =\cup_{n=0}^\infty \Ga^{(n)} (\X),
$$
where $\Ga^{(n)}(\X)$ denoted the set of all configurations with $n$ elements.
We will see that in the continuous combinatoric the spaces $\Ga(\X)$ and $\Ga_0(\X)$ 
will play very different roles. It is a specific moment related with transition to the continuum.
In this sense $\N$ is splitting  in these two spaces of configurations that makes corresponding
combinatorics  essentially more reach and sophisticated.

Configuration spaces present beautiful combinations of discrete and continuous properties. In particular,
in these spaces we have interesting differential geometry, differential operators and diffusion processes
etc., see e.g. \cite{AKR}.  From the other hand side, discreteness of an individual configuration makes possible to introduce 
proper analog of the difference calculus.

Note from the beginning, that the analog of the extension $\N\subset \R$ now naturally play the 
pair $\Ga(\X)\subset \mathcal{M}(\X)$ where we have in mind an imbedding of configurations
in the space of discrete  Radon measures on  $\X$ and, as a result, in the  space of 
all Radon mesures on $\X$ :

$$
\Ga(X) \ni  \ga \mapsto \ga(dx)= 
\sum_{y\in\ga} \delta_y (dx) \in \M(\X).
$$

Therefore, instead of pair
$$
\N\subset \R
$$
we have 
$$
\Ga(\X) \subset \mathcal{M}(\X).
$$
As a result, the transition to "continuous" variables in the 
considered situation leads to functions on $\mathcal{M}(\X)$.
In spatial combinatorics many  objects will be measure-valued.

Now  we will introduce an analog of  the generation function from
classical combinatorics. 
For a test function from the Schwarz space of test functions $\mathcal{D}(\X)$ 
$0\leq \xi \in  \mathcal{D}(\X)$ consider a function
$$
E_{\xi}  (\omega) = e^{<\ln (1+\xi),\omega>} \;\; \omega \in \mathcal{D}'(\X)
$$
that is a function on the space of Schwarz distributions. 
The power decomposition w.r.t. $\xi$ gives
$$
E_\xi(\omega)=\sum_{n=0}^\infty  \frac{1}{n!}  <\xi^{\otimes n}, (\omega)_n>.
$$
Generalized kernels $(\omega)_n \in \mathcal{D}'(\R^{nd})$ are called infinite dimensional falling factorials on
$\mathcal{D}'(\X)$. Define binomial coefficients (Newton polynomials)  on $ \mathcal{D}'(\X)$ as
$$
\binom \omega n= \frac{(\omega)_n}{n!}.
$$
Note that these objects are defined now on the very big space of distributions.
In particular cases we shall restrict them on the space of configuration
or Radon measures.

In particular, infinite dimensional Chu-Vandermond relations on configurations is
$$
(\ga_1 \cup \ga_2 )_n= \sum_{k=0}^n \binom nk (\ga_1)_k \otimes (\ga_2)_{n-k}.
$$

\begin{theorem}
For $\om \in \mathcal M (\X)$ 
$$
(\om)_0= 1
$$
$$
(\om)_1 =\om
$$
$$
(\om)_n(x_1,\dots,x_n)= \om(x_1) (\om(x_2)-\delta_{x_1}(x_2))
\dots (\om(x_n) -\delta_{x_1}(x_n) - \dots -  \delta_{x_{n-1}} (x_n)).
$$
\end{theorem}

 In the particular case $\om=\ga =\{x_i\; |\; \; i\in\N\}$
 
 $$
 (\ga)_n  = n! \binom \ga n = \sum_{\{i_1.\dots, i_n\} \subset \N}
 \delta_{x_1} \odot \dots \odot \delta_{x_n},
 $$ 
where $\delta_{x_1} \odot \dots \odot \delta_{x_n}$ denotes symmetric tensor
product.

We have 
$$
\Ga(\X) \ni \ga \mapsto \ga(dx) \in \M(\X).
$$
Due to our construction
$$
(\ga)_n \in \M(\R^{nd})
$$
is a symmetric Radon masure.
Therefore, we arrive in measure valued Newton polynomials.
The latter is the main consequence of  continuous combinatoric transition.

\section{Difference geometry for spatial combinatorics}

For any $x\in \ga$ define an elementary Markov death operator (death gradient)
$$
D^-_x F(\ga)= F(\ga\setminus x)-F(\ga)
$$
and the tangent space $T^-_{\ga} (\Ga)= L^2 (\R^d, \ga)$. Then for $\psi \in C_0(\R^d)$
$$
D^-_{\psi} F(\ga)=\sum_{x\in\ga} \psi(x) (F(\ga\setminus x)-F(\ga))
$$
is the directional (difference)  derivative.

Similarly, we define for $x\in \R^d$ 
$$
D^+_x F(\ga)= F(\ga \cup x)-F(\ga)
$$
and the tangent space $T^-_{\ga} (\Ga)= L^2(\R^d, dx)$. Then for $\varphi \in C_0(\R^d)$ 
$$
D^+_\varphi F(\ga)= \int_{\R^d} \varphi(x) (F(\ga \cup x)-F(\ga)) dx
$$
is another directional (difference) derivative.

For $\varphi \in C_0(\R^d)$ define a function
$$
E_\varphi (\ga) = \exp(<\ga, \log(1+\varphi)>),\;\;\ga \in \Ga.
$$
It is the generation function for the system on falling factorials (Newton polynomials) on $\Ga$:
$$
E_\varphi (\ga) = \sum_{n=0}^{\infty} \frac{1}{n!} <\varphi^{\otimes n}, (\ga)_n>.
$$
Then 
$$
D^+_\psi E_\varphi(\ga)=  <\varphi \psi> E_\varphi (\ga).
$$
An explicit formula for the falling factorials (as measures on $(\R^d)^n$) is
$$
(\ga)_n = \sum_{\{x_1,...,x_n\} \subset \ga  } \delta_{x_1} \odot 
\delta_{x_2}\odot \cdots \odot \delta_{x_n},
$$
where $\delta_{x_1} \odot 
\delta_{x_2}\odot \cdots \odot \delta_{x_n}$ denotes the symmetric tensor product of measures.

The action of difference derivatives on Newton monomials is given by
$$
D^+_{\psi} <\varphi^{(n)}, (\ga)_n>= n \int_{\R^d} \psi(x)<\varphi^{(n)}(x,\cdot), (\ga)_{n-1}(\cdot)>dx,
$$
$$
D^-_{\psi} <\varphi^{(n)}, (\ga)_n>= -n \sum_{x\in\ga} \psi(x) <\varphi^{(n)}(x,\cdot), (\ga \setminus x)_{n-1}(\cdot)>.
$$

\section{Stirling kernels}

We have polynomial equality 
$$
(\ga)_{n}= \sum_{k=1}^n s^n_k  \ga^{\otimes k},
$$
where
$$ 
s_k^n: \mathcal{D}'(\R^{kd}) \to \mathcal{D}'(\R^{nd})
$$
is a linear mapping.

On other side
$$
\ga^{\otimes n} =\sum_{k=0}^n S^n_k (\ga)_k,
$$
where
$$
S_k^n :  \mathcal{D}'(\R^{kd}) \to \mathcal{D}'(\R^{nd})
$$
is a linear mapping.

Kernels $s_k^n$ and $S^n_k$ we will call Stirling kernels of
first and second kind respectively. In the classical combinatorics Stirling 
coefficients play  a very important role.

For $f^{(n)}\in \mathcal{D}(\R^{nd})$
$$
<(\ga)_n, f^{(n)} > = 
$$
$$
\sum_{k=0}^n \frac{n!}{k!} 
< \ga^{\otimes k}(x_1,\dots,x_k),
\sum_{i_1 + \dots i_k =n} 
\frac{(-1)^{n+k}}{i_1\dots i_k} f^{(n)} (\underbrace{x_1,\dots,x_1}_{i_1 \rm{times}}, \dots, \underbrace{x_k,
\dots,x_k}_{i_k \rm{times}}) >.
$$

For the second kind kernels

$$
< \ga^{\otimes n}, f^{(n)}>=
$$
$$
\sum_{k=0}^{n} \frac{1}{k!} < (\ga)_k(x_1,\dots,x_k) , \sum_{i_1 + \dots i_k =n} 
\binom{n}{i_1\dots i_k} f^{(n)} (\underbrace{x_1,\dots,x_1}_{i_1 \rm{times}}, \dots, \underbrace{x_k,
\dots,x_k}_{i_k \rm{times}}) >.
$$

\section{Harmonic analysis on $\Ga(\X)$}

Functions $G:\Ga_0(\X) \to \R$ we call quasi-observables.
Note that $G$  restricted on $\Ga^{(n)} (\X)$ is given by a symmetric kernel
$G^{(n)}(x_1,\dots,x_n)$  and then
$$
G=(G^{(n)})_{n=0}^\infty.
$$

Functions $F: \Ga(\X) \to\R$ we call observables. For a quasi-observable $G$ define an operator

$$
(KG)(\ga) = \sum_{\eta\subset \ga,
|\eta|<\infty} G(\eta), \;\; \ga\in \Ga(\X)
$$
that is  an observable. To be well defined we need certain assumptions about
$G$ \cite{KoKu}.

For $G_1,G_2: \Ga_0(\X)\to \R $
define
$$
(G_1\star G_2)(\eta)= \sum_{\eta_1 \cup \eta_2\cup \eta_3 =\eta}
G_1(\eta_1\cup \eta_2) G_2(\eta_2\cup \eta_3) .
$$
Then
$$
K(G_1\star G_2)= KG_1  KG_2.
$$

Let $\mu\in \M^1 (\Ga(\X))$. 

$$K: Fun (\Ga_0) \to Fun (\Ga)$$,
$$
K^* : \M^1(\Ga) \to \M(\Ga_0).
$$
$$
K^* \mu =\rho,\;\; \rho= (\rho^{(n)} )_{n=0}^\infty.
$$
The measure $\rho$ is called correlation measure for $\mu$
(Fourier transform of $\mu$).

Assume absolute continuity
$$
d\rho^{(n)}(x_1,\dots ,x_n) = \frac{1}{n!} k^{(n)} (x_1,\dots, x_n) dx_1 \dots x_n.
$$
We call $ k^{(n)} (x_1,\dots, x_n),n\in \N$ correlation functions of the measure $\mu$.

Transition from measures to CFs is one of the main technical aspects of the analysis on 
CS in applications to dynamical problems.

Alternatively define the Bogoliubov functional
$$
B_\mu (\phi) =\int_{\Ga(\X)} e^{<\ga,\log(1+\phi>)} d\mu(\ga).
$$
Assuming  $B_\mu$ is holomorphic  in $\phi\in L^1 (\X)$
we obtain 

$$
B_\mu (\phi) =\sum_{n=0}^\infty \frac {1}{n!} \int k^{(n)}(x_1,\dots,x_n) \phi(x_1) \dots \phi(x_n) dx_1\dots dx_n.
$$

\section{From spatial to classical combinatorics}

Having developed combinatorial structures in the continuum, we 
may consider  the inverse direction. Namely, how looks like our infinite-dimensional
objects in the one dimensional reduction. Surprisingly, it may give some new structures
even in this classical case.

Let $a,b: \N\to \R$. 
Define a convolution
$$
(a\star  b)(n) =\sum_{j+k+l=n} a(j+k) b(k+l) .
$$
As before
$$
(Ka)(n)  = \sum_{k+0}^n \binom nk a(k).
$$
Then 
$$
K(a\star b) = Ka \cdot Kb.
$$

Introduce coherent states
$$
e_\la(\cdot) : \N\to \C,\; e_\la(n)= \la^n,\;\;\la\in \C.
$$
$$
(Ke_\la)(n)= (1+\la)^n.
$$

The configuration space $\Ga(\X)$ is the space of microscopic states in the 
classical statistical physics of continuous systems.
A measure $\mu\in\M^1(\Ga(\X))$ is a macroscopic state of a continuous system in 
the statistical physics.
Coming back we can interpret (a bit naively) a measure $\mu\in \M^1(\N)$ as  a state of $0D$ system.

For example, the Poisson measure for $\sigma >0$
is defined as 
$$
\pi_\sigma (n)= e^{- \sigma} \frac{\sigma^n}{n!}.
$$

Several characteristics we can incorporate in such a case  from the analysis on $\Ga(\X)$.
Introduce the Bogoliubov functional:

$$
B(\la)= \int_{\R_+}  (1+\la)^x d\mu(x).
$$

$$
(1+\la)^x= \sum_{n=0}^\infty  \frac{\la^n}{n!} (x)_n.
$$

\begin{theorem}

Let $\mu \in \M^1(\R_+)$.  Then $\mu(\N)=1$ iff
$B(\la)$  has a holomorphic  extension.

\end{theorem}
Similarly we can define correlation measures

$$
\int_\N (Ka) (x) d\mu(x) =\int_\N a(x) d\rho_\mu (x).
$$

$$
\rho_\mu (n) = \frac{1}{n!} \int_\N (x)_n d\mu(x) =
$$
$$
\sum_{m=n}^\infty \binom mn \mu(m).
$$

\vspace{5mm}
\section{Acknowledgment}

The financial support by the Ministry 
for Science and Education of Ukraine
through Project 0119U002583
is gratefully acknowledged.

\end{document}